\documentclass[a4paper,notitlepage, 12pt,reqno]{article}

  \usepackage{mathtext}
       \usepackage{graphicx}
 \usepackage[cp1251]{inputenc}
  \usepackage[T2A]{fontenc}
 \usepackage[english]{babel}
 \usepackage[all,arc,poly,2cell,curve,arrow,tips]{xy}

  \usepackage{enumerate, eucal, amsthm,amsmath, amssymb}

  \allowdisplaybreaks[4]

 \theoremstyle{definition}
 \newtheorem{defn}{Definition}

 \theoremstyle{plain}
 \newtheorem{thm}{Theorem}

 \newtheorem{prop}{Proposion}

\newtheorem{cor}{Corrolary}

 \newtheorem{lem}{Lemma}

 \theoremstyle{remark}

 \renewcommand{\abstractname}{}

  \newcounter{ab}
\title{
Noncommutative pfaffians and classification of states of
five-dimensional quasi-spin.}
\author{Dmitry Artamonov,  Valentina Goloubeva}

\begin{document}

\maketitle

\renewcommand{\abstractname}{}

\begin{abstract} Noncommutative pfaffians associated with an orthogonal  algebra $\mathfrak{o}_N$ are some
 special elements of the universal enveloping algebra $U(\mathfrak{o}_N)$.
Using  pfaffians we construct the fourth quantum number  which
together with the naturally defined three quantum numbers  allow to
classify  the states of a five-dimensional quasi-spin. The pfaffians
are treated as creation operators for the new quantum number.
\end{abstract}

\maketitle{}

\section{The quasi-spin algebra}

In many  problems of quantum numbers  there appears naturally an
algebra of observables whose elements are quadratic combinations of
fermion creation and annihilation operators. Such  an  algebra is
called a generalized fermion algebra. Some its particular cases are
called the quasi-spin algebra \cite{Heb}.

For example,  the three-dimensional  quasi-spin algebra appears
naturally when one considers particles of the same type. The three
dimensional quasi-spin  is constructed as follows. One takes the
system of noninteracting fermions with the fixed angular momentum
$j$. Denote by $a_m^+$, $a_m$,  $m=-j,...,j$, the creation and
annihilation operators of particles with the projection $m$ of the
angular momentum to the axes
 $z$.

 Consider the operators

 $$s_k^+=a_k^+a_{-k}^+,$$

 $$s_{k}^-=a_{-k}a_k,$$

 $$s^0_k=\frac{1}{2}(a_k^+a_k-a_{-k}a_{-k}^+).$$

Let us construct a three-dimensional algebra of quasi-spin,
isomorphic to $\mathfrak{o}_3=\mathfrak{sl}_2$, see
\cite{Lip},\cite{Bid}. It is spanned by elements

$$S_+=\sum_k s_k^+,$$

$$S_-=\sum_k s_k^-,$$

$$S_0=\sum_k s_k^0.$$

For the classification of states of the three-dimensional quasi-spin
the quantum numbers seniority and the number of particles are used.
Note that although the three-dimensional quasi-spin algebra is
isomorphic to the spin algebra the two quantum  numbers  mentioned
above are different from the quantum numbers used for the
classification of the states of the spin algebra.


The three-dimensional quasi-spin algebra appears  in the theory
superconductivity, in description of pair correlations between
nucleons, \cite{Bid}.

Suppose we have particles of two types, protons and neutrons.
 Then the five-dimensional quasi-spin algebra naturally appears.
 Each creation and annihilation operators has two indices.
  Let $a_{pm}^+$, $a_{nm}^+$ be  creation operators of protons and
  neutrons with the projection $m$ of angular momentum to the $z$-axes, and $a_{pm}$, $a_{nm}$ the corresponding annihilation operators.
   The commutation relations between these operators are the following (the index $\tau$ is either  $p$ or  $n$):

$$a_{\tau m}a_{\tau'm'}+a_{\tau' m'}a_{\tau m}=0$$

$$a_{\tau m}^+a_{\tau'm'}^++a_{\tau' m'}^+a_{\tau m}^+=0$$

$$a_{\tau m}a^+_{\tau'm'}+a_{\tau' m'}^+a_{\tau m}=\delta_{\tau,\tau'}\delta_{m,m'}$$

In this case one can naturally define the five-dimensional
quasi-spin algebra \cite{Helme} (see also \cite{Lip}), spanned by
elements

$$\tau_+=\sum_m a_{pm}^+a_{nm},\,\,\,\, \tau_0=\frac{1}{2}\sum_m(a_{p m}^+a_{p m}-a_{p m}^+a_{n m}),\,\,\,\, \tau_{-}=\sum_m a_{nm}^+a_{p m}$$

$$N=\frac{1}{2}\sum_m(a_{pm}^+a_{pm}+a_{nm}^+a_{nm})-\frac{2j+1}{2}$$

$$A(1)=\sum_{m>0}(-1)^{j-m}a_{pm}^+a_{p-m}^+,\,\,\,
A(0)=\frac{1}{\sqrt{2}}\sum_{m>0}(-1)^{j-m}(a^+_{pm}a^+_{n-m}+a^+_{n-m}a^+_{p-m}),$$$$
A(-1)=\sum_{m>0}(-1)^{j-m}a^+_{nm}a^+_{n-m}$$

$$B(1)=\sum_{m>0}(-1)^{j-m}a_{n-m}a_{pm},\,\,\,
B(0)=\frac{1}{\sqrt{2}}\sum_{m>0}(-1)^{j-m}(a_{n-m}a^+_{pm}+a^+_{p-m}a^+_{nm}),$$$$
B(-1)=\sum_{m>0}(-1)^{j-m}a_{n-m}a_{nm} $$

 These operators form a Lie algebra isomorphic to $\mathfrak{o}_5$. Note that the subalgebra $\mathfrak{o}_3$
 spanned by
$\tau_+,\tau_0,\tau_-$ is the isospin subalgebra and $N$ is the
number of particles plus a constant. But below for simplicity of
terminology we call $N$ itself the number of particles. The reasons
for the introduction of this five-dimensional quasi-spin are
discussed in \cite{He1}.

  One can obtain an  explicit isomorphism of the quasi-spin algebra with the algebra
    $\mathfrak{o}_5$  in the split realization (see Sec. \ref{or}) in the following way. First one must establish a relation between the
    written above
      elements $\tau_i, A(j),B(l)$ and the Chevalley
    base of
 $\mathfrak{o}_5$. Second one must establish a relation between  canonical generators of  $\mathfrak{o}_5$  and the Chevalley
    base of
 $\mathfrak{o}_5$.  The  commutation relation between   elements $\tau_i, A(j),B(l)$  that are necessary for the above procedure
 are given in  \cite{Heb}. As the result we have

$$F_{0-1}=\frac{\tau_+}{\sqrt{2}}\,\,\, F_{-1-2}=A(-1),\,\,\, F_{0-2}=A(0),\,\,\, F_{2-1}=-A(1),\,\,\, F_{-1-1}=-\tau_0,$$

$$F_{-10}=\frac{\tau_-}{\sqrt{2}}\,\,\, F_{-2-1}=B(-1),\,\,\, F_{-20}=B(0),\,\,\, F_{-12}=-B(1),\,\,\, F_{-2-2}=-N,$$

where $F_{ij}$ are defined in Sec. \ref{or}.


Now let us proceed to the problem of classification of states  of a
finite-dimensional representation of the five-dimensional
 the  quasi-spin algebra.

 There are three natural quantum
 numbers indexing the base vectors of a representation. They are the
 number of particles $N$, the isospin $T$ and the projection of the isospin $\tau_0$.

These numbers are not enough to classify the states, there are
linear independent states with the same collections of these three
numbers.

Thus, in order to classify the states it is necessary to construct
at least one new quantum number $k$.

First construction of such a number in \cite{Goli},\cite{FS} were of
the following type.

 At first one fixes the set of states, to which
the zero value of the fourth quantum number is assigned. Then the
creation operator for the new quantum number is introduced.  Then
one defines the fourth quantum number inductively in the following
manner. If a state is obtained from a state with the fourth quantum
number $k$ by application of the creation operator for the fourth
quantum number, then we assign to it the fourth quantum number
$k+1$.

However, it turns out that the creation operator of the new quantum
number changes other quantum numbers in a not clear way.

Later other constructions of the fourth quantum number were given in
\cite{He1},\cite{He2},\cite{Sz},\cite{A},\cite{ST}. There exist
other solutions of the same classification problem  that use the
technique of the vector coherent states
\cite{Rore1},\cite{Rore2},\cite{Rore3}.

Note that the problem of classification of states of the
five-dimensional quasi-spin appears not only in the shell model of
nuclear structure described above but also in  of classification of
states of two-dimensional oscillator
\cite{Goli},\cite{He1},\cite{He2}, in the description of states of
the Bohr-Mottelson model
  \cite{Mod11} in the model of interacting bosons \cite{Mod12} and in
  others (see references in
  \cite{Rore3}).

In the present paper another solution of the problem of construction
of the fourth quantum number is given.  The quantum number
constructed in the paper has several advantages over quantum numbers
constructed  earlier. As in \cite{Goli},\cite{FS}  the fourth
quantum number is introduces using the creation operator for this
quantum number. This allows to give a simple physical interpretation
for the new quantum number. But  in contrast to
\cite{Goli},\cite{FS} our  creation operator changes other quantum
numbers in a very simple way.

We use as creation operator the noncommutative pfafians associated
with the algebra $\mathfrak{o}_5$.

The paper consists of two parts. In the first part we give a
construction of the fourth quantum number using noncommutative
pfaffians. In this construction and in the proof that our quantum
number does solve the problem of classification of states a
technical theorem is used. The second part of the paper, named the
appendix, is devoted to the proof of this theorem.

\section{The orthogonal algebra and nocommutative pfaffians}

\label{or}

To present our construction of the fourth quantum number we must
introduce noncommutative pfaffians and  explain, what is the
Gelfand-Tsetlin-Molev base of a $\mathfrak{o}_5$-representation.

Let  $\Phi=(\Phi_{ij})$, $i,j=1,...,2n$  be a skew-symmetric
 $2n\times 2n$-matrix, whose  matrix entries belong to
  a noncommutative ring.

\begin{defn}The noncommutative pfaffian
 of $\Phi$ is defined by the formula

$$Pf \Phi=\frac{1}{n!2^n}\sum_{\sigma\in
S_{2n}}(-1)^{\sigma}\Phi_{\sigma(1)\sigma(2)}...\Phi_{\sigma(2n-1)\sigma(2n)},$$

were $\sigma$ is a permutation of the set  $\{1,...,2n\}$.
\end{defn}

In the paper a split realization of $\mathfrak{o}_N$ is used. To
formulate it we use the following indexation of rows and columns of
matrices from  $\mathfrak{o}_N$.  When $N$ is odd the indices $i,j$
of rows and columns belong to the set $\{-n,...,-1,0,1,...,n\}$,
where $n=\frac{N-1}{2}$. When $N$ is even the indices $i,j$ belong
to the set $\{-n,...,-1,1,...,n\}$, where $n=\frac{N}{2}$. Shortly
in both cases this set of indices is denoted in the paper as
$\{-n,...,n\}$.

Then the algebra $\mathfrak{o}_{N}$ is  defined as the Lie algebra
spanned  by matrices $F_{ij}=E_{ij}-E_{-j-i}$, where $E_{ij}$ are
matrix units.

One can prove that elements $F_{-n-n},...,F_{-1-1}$ form a base in
the Cartan subalgebra, and  elements $F_{ij},j<-i$ are root
elements.

The precise   correspondence is the following. Let $e_i$ be the
element  $F_{ii}^*$ in the dual space to the Cartan subalgebra. Put
$e_{-r}:=-e_{r}$ and $e_0=0$.  Then the element  $F_{ij}$
corresponds to the root $e_i-e_{j}$.

The commutation relations between the generators are the following

$$[F_{ij},F_{kl}]=\delta_{kj}F_{il}-\delta_{il}F_{kj}-\delta_{-ki}F_{-jl}
+\delta_{-lj}F_{k-i}.$$

In the paper the following noncommutative pfaffians are  considered.

\begin{defn} Let $F$ be the matrix $F=(F_{ij})$.
For every subset $I\subset\{-n,...,n\}$ which consists of an even
number $k$ of elements define a submatrix $F_I$   by the formulae
$F_{I}=(F_{ij})_{-i,j\in I}$. Put

$$PfF_I=Pf(F_{-ij})_{-i,j\in I}.$$
\end{defn}

In \cite{Molev}  the author  in terms of these pfaffians defines
some special elements of $U(\mathfrak{o}_N)$ called the Capelli
elements
 $C_k=\sum_{I\subset \{-n,...,n\},|I|=k}PfF_IPfF_{-I}$,
$k=2,4,...,[\frac{N}{2}]$.  It is proved that the elements $C_k$
belong to the center of $U(\mathfrak{o}_N)$.

In the case  $N=2n+1$ put
$$PfF_{\widehat{-n}}:=PfF_{\{-n+1,...,n\}}\text{ and }
PfF_{\widehat{n}}:=PfF_{\{-n,...,n-1\}}$$




If one identifies  $\mathfrak{o}_5$ with the quasi-spin algebra then
the pfaffians $PfF_{\widehat{2}}$, $PfF_{\widehat{-2}}$ are written
as follows
$$PfF_{\widehat{2}}=A(-1)\star\frac{\tau_+}{\sqrt{2}}+A(0)\star\tau_0+A(1)\star\frac{\tau_-}{\sqrt{2}},$$

$$PfF_{\widehat{-2}}=B(-1)\star\frac{\tau_-}{\sqrt{2}}+B(0)\star\tau_0+B(1)\star\frac{\tau_+}{\sqrt{2}},$$

where $a\star b=\frac{1}{2}(ab+ba)$.

\section{The Gelfand-Tsetlin-Molev base of a $\mathfrak{o}_{5}$-representation}
\label{dpb}

For the construction of the fourth quantum number we use the
Gelfand-Tsetlin-Molev base of a $\mathfrak{o}_{5}$-representation.
In this section we give its definition and describe the action of
the pfaffian  $PfF_{\widehat{2}}$ in this  base.

 The Gelfand-Tsetlin-Molev base of a  $\mathfrak{o}_{2n+1}$-representation is a base of a Gelfand-Tsetlin type, whose construction is  based
   on restrictions
 $\mathfrak{o}_{2n+1}\downarrow \mathfrak{o}_{2n-1}$, in contrast to the classical Gelfand-Tsetlin base,
 whose construction is based on restrictions
  $\mathfrak{o}_N\downarrow \mathfrak{o}_{N-1}$.

One can give the following formal definition. The
Gelfand-Tsetlin-Molev base is a  weight base of a
$\mathfrak{o}_{2n+1}$-representation and for different
  $n$  the procedures of constructions of  such bases must be coherent  in the following sense. A base of a  $\mathfrak{o}_{2n+1}$-representation
  must be a union of bases in  $\mathfrak{o}_{2n-1}$-representation,
  into which a
  $\mathfrak{o}_{2n+1}$-representation splits when one restricts $\mathfrak{o}_{2n+1}\downarrow \mathfrak{o}_{2n-1}$.

Note that in the case $\mathfrak{sp}_{2n}$ an analogous base was
firstly constructed by Zhelobenko \cite{zb}.

Describe a base of a  $\mathfrak{o}_5$-representation. All notations
are taken from \cite{Molev}. In construction of such a base only one
restriction $\mathfrak{o}_5\downarrow \mathfrak{o}_3$ appears. Hence
base vectors are weight vectors for the algebra $\mathfrak{o}_5$ and
every base vector is contained in a $\mathfrak{o}_3$-representation
that appears, when one restricts $\mathfrak{o}_5\downarrow
\mathfrak{o}_3$.

The base vectors of a $\mathfrak{o}_5$-representation with the
highest weight $(\lambda_1,\lambda_2)$,
$$0\geq \lambda_1\geq \lambda_2$$ are indexed by tableaus $\Lambda$
of type

$\lambda_1,\lambda_2$

$\sigma,\lambda'_{21},\lambda'_{22}$

$\lambda_{11}$

$\sigma_1,\lambda'_{11}$

Let us give an interpretation of these numbers and give inequalities
for them.

\begin{enumerate}
\item  The number $\lambda_{11}$ is a weight  of a  $\mathfrak{o}_3$-representation $\lambda_{11}$ that contains the base vector. This number must satisfy the inequalities
 $$0\geq \lambda_1\geq \lambda_{11}\geq \lambda_2.$$
\item The numbers $\sigma_1$ and $\lambda'_{21}$ give an index of a base vector in a $ \mathfrak{o}_3$-representation with the
highest weight $\lambda_{11}$ that contains the considered base
vector. The constraints on these numbers are
$$0\geq \lambda'_{21}\geq\lambda_{21}$$ and $\sigma_1=0,1$.
\item The numbers   $\sigma$, $\lambda'_{21},\lambda'_{22}$ define an element of a base
in the space of $\mathfrak{o}_3$-highest vectors with the
$\mathfrak{o}_3$-weight $\lambda_{11}^{+}$.  These numbers satisfy
the constraints
   $$0\geq \lambda'_{21}\geq \lambda_1\geq \lambda'_{22}\geq \lambda_2,$$ $$0\geq \lambda'_{21}\geq \lambda_{11}\geq
   \lambda'_{22},$$ $$\sigma=0,1; \text{ if
    }\lambda'_{21}=0,\text{ then }\sigma=0.$$
\end{enumerate}

All numbers $\lambda$ are simultaneously integer or half integer.

Let us establish a relation between quantum numbers $N,T,\tau_0$ and
elements of the tableau $\Lambda$.

The number $T$  is a highest weight of a
$\mathfrak{o}_3$-representation that contains the base vector. Thus
$$T=\lambda_{11}.$$

The number $\tau_0$ is the index of the base vector in it
$\mathfrak{o}_3$-representation that contains it. In the
Gelfand-Tsetlin-Molev base an element of a  $
\mathfrak{o}_3$-representation with the highest weight $0\geq
\lambda_{11}$ is indexed by numbers $\sigma_1=0,1$
  and $\lambda'_{11}$, such that  $0\geq \lambda'_{11} \geq \lambda_{11}$.
   Traditionally the elements of a weight base are indexed by one number    $\tau_0$,  such that $\lambda_{11}\leq \tau_0 \leq -\lambda_{11}$.
    The relation between these approaches is the following
     $$\tau_0=(-1)^{\sigma_1}\lambda'_{11}.$$

The  number $N=F_{22}$ is expressed
as follows \cite{Molev}

$$N=\sigma+2(\lambda'_{2,1}+\lambda'_{2,2})-(\lambda_1+\lambda_2)-\lambda_{1,1}$$

From the mathematical point of view the problem  of construction of
the fourth quantum number is a problem of  construction of an
indexation of the vectors of a base space of
$\mathfrak{o}_3$-highest vectors with the $\mathfrak{o}_3$-weight
$T$  (this space is called the multiplicity space) using the number
$N$ and some new quantum number.

Using $\sigma,\lambda'_{2,1},\lambda'_{2,2}$ one can solve the
problem of construction of the fourth quantum number.  But then the
physical interpretation of the obtained quantum number is difficult.

Introduce notation $\rho_i=i-\frac{1}{2}$ for $i>0$, and
$\rho_{-i}=-\rho_{i}$.

Define $$\gamma_i=\lambda'_{2i}+\rho_i+\frac{1}{2}.$$

Also put $$\overline{\sigma}=\sigma+1 \,mod\, 2.$$

\begin{thm}
\label{maint1} On the vector $\xi_{\Lambda}$ the pfaffian
$PfF_{\widehat{2}}$ acts as follows.

Let $\Lambda=(\lambda,\sigma,\lambda',\Lambda')$, where $\lambda$ is
the first row of   $\Lambda$, $\sigma,\lambda'$ is the second row of
$\Lambda$ and $\Lambda'$ is the remaining part of $\Lambda$.

Let $\xi_{\sigma,\lambda',\Lambda'}$ be the corresponding base
vector.

If $\sigma=0$, then
$$PfF_{\widehat{2}}\xi_{\sigma,\lambda',\Lambda'}=\xi_{\bar{\sigma},\lambda,\Lambda'}.$$


If $\sigma=1$, then
$$PfF_{\widehat{2}}\xi_{\sigma,\lambda',\Lambda'}=(-1)^n\sum_{j=1}^2\Pi_{t=1,t\neq
j}^2\frac{-\gamma_t^2}{\gamma_j^2-\gamma_t^2}\xi_{\bar{\sigma},\lambda'+\delta_j,\Lambda'},$$

Where $\lambda'+\delta_j$ is the row $\lambda'$ with $1$ added to
the $j$-th component.

\end{thm}

The  general form of the theorem \ref{maint1} (the theorem
\ref{maint}) is proved in the appendix.

\section{Construction of an additional quantum number using
pfaffian}

\label{number}

Let us show how to construct the fourth quantum number using
pfaffians.


Denote as $|l_1,...,l_p>$ a vector corresponding to quantum numbers
$l_1,...,l_p$.

The main result of this section is the following.

\begin{thm}\label{teork}
There exists an indexation  of base weight vectors of a
representation of $\mathfrak{o}_5$ by numbers $T,\tau_0,N,k$, where
the additional number  $k$ is a nonnegative integer. States which
have different collections of indices are independent.

If $N<0$, then $PfF_{\widehat{2}}$ maps $|T,\tau_0,N,k>$ to
$|T,\tau_0,N+1,k+1>$.

If $N>0$, then $PfF_{\widehat{-2}}$ maps $|T,\tau_0,N,k>$ to
$|T,\tau_0,N-1,k+1>$.
\end{thm}

All the rest part of the section is devoted to the construction of
the quantum number $k$ with the prescribed properties and to the
proof that the numbers $T,\tau_0,N,k$  are sufficient for the
classification the base vectors of a
$\mathfrak{o}_5$-representation.

\proof

First of all we construct the additional number for the state for
with
 $N\leq 0$  and prove that the new quantum number solves the problem of classification.
 Then  for states with $N>0$ the quantum number is constructed using the reflection of the algebra $\mathfrak{o}_5$.

The construction of the quantum number and the proof that it solves
the classification problem is done by induction by  $N$, the numbers
 $T,\tau_0$ are suggested to be fixed.

The base vectors of the Gelfand-Tsetlin-Molev base of a
$\mathfrak{o}_5$-representation are encoded by tableaus, which
include actually the numbers  $\lambda''_{11},\sigma_1$ related to
$T,\tau_0$, and also the numbers
$\sigma,\lambda'_{2,1},\lambda'_{2,2}$. The restriction on
$\sigma,\lambda'_{2,1},\lambda'_{2,2}$ are the following
$$0\geq \lambda'_{2,1}\geq \lambda_1\geq \lambda'_{2,2}\geq
\lambda_2,$$
$$\lambda'_{2,1}\geq \lambda_{1,1}\geq \lambda'_{2,2},$$
$$\sigma=0,1.$$

In other words the point $(\lambda'_{2,1},\lambda'_{2,2})$ belongs
to the rectangle on the figure \ref{fig1}.

\begin{figure}
	\label{fig1}
	$$
	{\tiny
		\xymatrix @ = 0.5pc @*[r]{
			\ar@{->}[dddddd]   \ar@{--}[rr]^{\max(\lambda_{11},\lambda_1)}    &&  \bullet \ar@{-}[rrrrr] \ar@{-}[dddd] &\bullet \ar@{-}[dddd] &\bullet\ar@{-}[dddd]& \bullet\ar@{-}[dddd] &\bullet\ar@{-}[dddd] & \bullet\ar@{-}[dddd] & \\
			&&  \bullet \ar@{-}[rrrrr]   &\bullet   &\bullet   & \bullet  & \bullet & \bullet  &    \\
			&&  \bullet \ar@{-}[rrrrr]   &\bullet   &\bullet   & \bullet  & \bullet & \bullet &  \\
			&& \bullet \ar@{-}[rrrrr] &\bullet   &\bullet   & \bullet  & \bullet & \bullet  &   \\
			\ar@{->}[rrrrrrrr]&   &\bullet   &\bullet   & \bullet  & \bullet & \bullet &  \bullet &
			\lambda'_{22} \\
			&_{\scriptsize {\min(\lambda_{11},\lambda_1)} } &&&&& &_{\lambda_2} & &  \\
			\lambda'_{21}  &  &&&&&&&&&
		}
	}
	$$
\end{figure}

\begin{figure}
\tiny{
	$$\xymatrix @ = 0.5pc @*[r]{
		\sigma=0\\
		&&  \bullet \ar@{-}[rrrrr] \ar@{-}[ddd] &\bullet \ar@{-}[ddd]  &\bullet\ar@{-}[ddd]&
		\bullet\ar@{-}[ddd] &\bullet\ar@{-}[ddd] & \bullet \ar@{-}[ddd] & \\
		&&  \bullet \ar@{-}[rrrrr]   &\bullet   &\bullet   & \bullet  & \bullet & \bullet  &    \\
		&&  \bullet \ar@{-}[rrrrr]   &\bullet   &\bullet   & \bullet  & \bullet & \bullet &  \\
		&& \bullet \ar@{-}[rrrrr] &\bullet   &\bullet   & \bullet  & \bullet & \bullet  &   \\
	}
	\;
	\xymatrix @ = 0.5pc @*[r]{\sigma=1\\
		&&  \bullet \ar@{-}[rrrrr] \ar@{-}[ddd] &\bullet \ar@{-}[ddd]  &\bullet\ar@{-}[ddd]&
		\bullet\ar@{-}[ddd] &\bullet\ar@{-}[ddd] & \bullet \ar@{-}[ddd] & \\
		&&  \bullet \ar@{-}[rrrrr]   &\bullet   &\bullet   & \bullet  & \bullet & \bullet  &    \\
		&&  \bullet \ar@{-}[rrrrr]   &\bullet   &\bullet   & \bullet  & \bullet & \bullet &  \\
		&& \circ  \ar@{-}[rrrrr] &\circ   &\circ   &\circ   & \circ  & \circ   &   \\
	}$$ }
\caption{The cases $\sigma=0$ and $\sigma=1$.}
\label{fig2}
\end{figure}
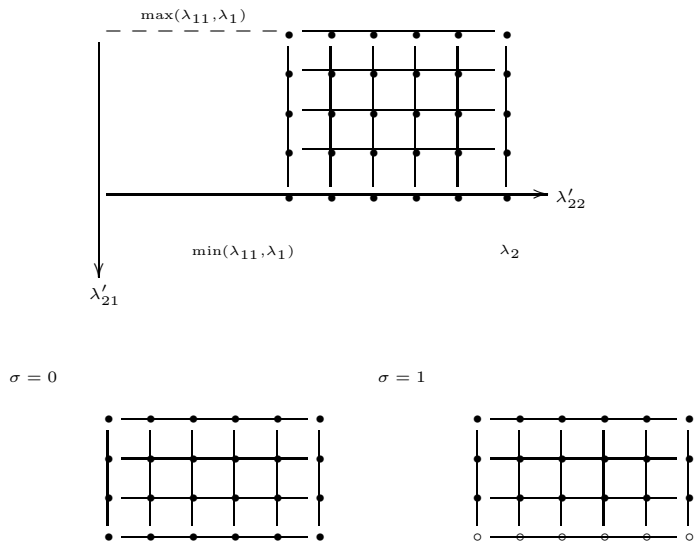

If $\sigma=0$ then the points on the lower edge are included, and if
$\sigma=1$ then they are not included, see figure \ref{fig2}.

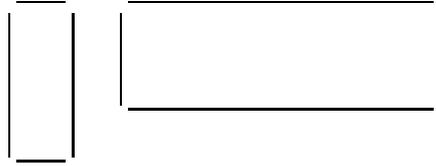
\begin{figure}
$$
\xymatrix @ = 0.5pc @*[r]{
	\ar@{-}[rr]  \ar@{-}[ddddd]   &     &   \ar@{-}[ddddd]   \\
	&&      \\
	&     &    \\
	&     &  \\
	&     &  \\
	\ar@{-}[rr] &&
}   \quad \xymatrix@R=0.5cm{
	\ar@{-}[rrrr]  \ar@{-}[dd]   &   &&   &   \ar@{-}[dd]   \\
	&     & && \\
	\ar@{-}[rrrr] && &&
}
$$
\caption{A narrow and a wide rectangle.}\label{fig3}
\end{figure}

The rectangle can have one of two forms, see figure \ref{fig3}.

In the first case we say that the rectangle is narrow and in the
second case we say that the rectangle is wide.

Write the formula for $N=F_{22}$ in the Gelfant-Tsetlin-Molev base
\cite{Molev}

$$N=\sigma+2(\lambda'_{2,1}+\lambda'_{2,2})-(\lambda_1+\lambda_2)-\lambda_{21}$$

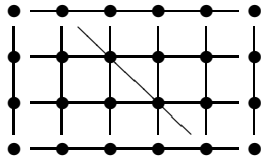
\begin{figure}
 \centering
 \centering
$$
\xymatrix @ = 0.5pc @*[r]{ N=const\\
	&&  \bullet \ar@{-}[rrrrr] \ar@{-}[ddd] &\bullet \ar@{-}[ddd] \ar@{-}[rrrddd]  &\bullet\ar@{-}[ddd]&
	\bullet\ar@{-}[ddd] &\bullet\ar@{-}[ddd] & \bullet \ar@{-}[ddd] & \\
	&&  \bullet \ar@{-}[rrrrr]   &\bullet   &\bullet   & \bullet  & \bullet & \bullet  &    \\
	&&  \bullet \ar@{-}[rrrrr]   &\bullet   &\bullet   & \bullet  & \bullet & \bullet &  \\
	&& \bullet \ar@{-}[rrrrr] &\bullet   &\bullet   & \bullet  & \bullet & \bullet  &   \\
}
$$
\caption{Points corresponding to the fixed $N$}\label{fig4}
\end{figure}

The point with fixed $N$ belong to the line
$\lambda'_{21}+\lambda'_{22}=const$, see figure \ref{fig4}.

\begin{figure}
\centering
$$
\xymatrix @ = 0.5pc @*[r]{  N=0\\
	\ar@{-}[rr]  \ar@{-}[ddddd]   &     &   \ar@{-}[ddddd]   \\
	\ar@{-}[rrddd]\\   &&      \\
	&     &    \\
	&     &  \\
	\ar@{-}[rr] &&
} \quad \xymatrix@R=0.5cm{ N=0\\
	\ar@{-}[rrrr]  \ar@{-}[dd]   & \ar@{-}[rrdd]   &&   &   \ar@{-}[dd]   \\
	&     & && \\
	\ar@{-}[rrrr] && &&
}
$$
\caption{Line $N=0$}\label{figu5}
\end{figure}
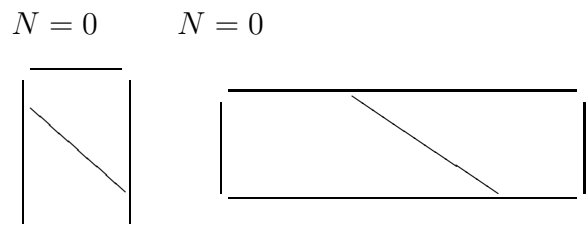

The line $N=0$ cuts the rectangle into two equal parts, see figure
\ref{figu5}.

Now we can begin the construction of the quantum number  $k$ for the
vectors with $N\leq 0$.

Remind that the construction is inductive by $N$ with fixed numbers
$T,\tau_0$.

{\bf  Base of induction.} Note that the state with $N=N_{min}$ and
fixed  $T,\tau_0$ is unique. This follows from the fact that to the
minimal value of $N$ there corresponds the upper right angle of the
rectangle. In other words for the minimal $N$  one has
$\lambda'_{2,1}=max\{\lambda_{1},\lambda_{1,1}\}$,
$\lambda'_{2,2}=\lambda_2$, $\sigma=0$.

Thus all numbers in the tableau $\Lambda$ are uniquely defined,
hence the states with $N=N_{min}$, and fixed $T,\tau_0$ is unique.

Assign to this state the zero number $k$. Since for $N=N_{min}$
there exist only one independent state for $N=N_{min}$ all states
are classified by four quantum numbers by obvious reasons.

{\bf  Preliminary discussion of the induction process. } Suppose
that for $N=N^*$ the states are classified by the number  $k$.
Consider the case $N=N^*+1$.

A base in the space of states with fixed  $T,\tau_0, N$  is  formed
by those vectors of the Gelfand-Tsetlin-Molev base, for which the
second row of the tableau $\Lambda$ has the form

$$(\sigma,K-t,t),$$ where $K$ and $N^*$ are related by $N^*=\sigma+2K-(\lambda_1+\lambda_2)-T$.

Denote such vectors $\xi^{N^*}_t$.

Note that for fixed $N$  one has for all states either  $\sigma=0$,
or $\sigma=1$.

The index $t$ runs only the values for which the point
$(\lambda_{2,1}',\lambda_{2,2}')=(K-t,t)$ belongs to the rectangle.


The action of the pfaffian $PfF_{\widehat{2}}$   on the second row
of the tableau $\Lambda$ is given  (see Theorem \ref{maint1}) by
formulae:

\begin{equation}\label{pkc1}(0,\lambda'_{2,1},\lambda'_{2,2}) \mapsto
(1,\lambda'_{2,1},\lambda'_{2,2})\end{equation}
\begin{equation}\label{pkc2}(1,\lambda'_{2,1},\lambda'_{2,2})\mapsto
c_1(0,\lambda'_{2,1}+1,\lambda'_{2,2})+c_2(0,\lambda'_{2,1},\lambda'_{2,2}+1),\end{equation}

where $c_1=\frac{-\gamma_2^2}{(\gamma_1^2-\gamma_2^2)}$,
$c_2=\frac{-\gamma_1^2}{(\gamma_2^2-\gamma_2^2)}$,
$\gamma_1=\lambda'_{2,1}$,
 $\gamma_2=\lambda'_{2,2}-1$.

If a row contained in these formulas does not satisfy the
constraints on the second row of a Gelfand-Tsetlin-Molev tableau
(that is the corresponding point does not belong to the rectangle),
then it must be replaced to zero.

 Since $\lambda'_{2,2}$ is non-positive, then $c_2\neq 0$, and
 $c_1=0$ if and only if $\lambda'_{2,1}=0$.
 But then  since $\lambda'_{2,1}+1>0$ the vector $(0,\lambda'_{2,1}+1,\lambda'_{2,2})$ which
 is multiplied by $c_1$ is zero. That is why we can  put in this case $c_1=1$ in the
formula \ref{pkc2} and suggest that always $c_1\neq 0$.

To the formulas \ref{pkc1} and \ref{pkc2} for the action of the
pfaffian $PfF_{\widehat{2}}$ the following figures correspond. To
the formula \ref{pkc1} there corresponds the figure \ref{fig6}.

\begin{figure}
\centering
$$
\xymatrix @ = 0.5pc @*[r]{
	\ar@{-}[rrrrr] &     &\ar@{-}[ddd] \ar@{-}[dddrrr]   &\ar@{-}[ddd]\ar@{-}[ddrr]&   \ar@{-}[d]   &  \ar@{-}[ddd] &\\
	&&   \ar@{-}[r]  &  \bullet     & \bullet    \ar@{=>}[l]\ar@{=>}[d] \ar@{-}[r] & &  \\
	&     &  \ar@{-}[rrr]&   & \bullet \ar@{-}[d]   &  & N=N^*,\; \sigma=1 \\
	&     &  \ar@{-}[rrr]&   &      && N=N^*+1,\; \sigma=0 }
$$
\caption{The action of the pfaffian $PfF_{\widehat{2}}$ in the case
$\sigma=1$}\label{fig6}
\end{figure}

\begin{figure}
\centering
$$N=N^*,\; \sigma=0
\xymatrix @ = 0.5pc @*[r]{
	\ar@{-}[rrrrr] &     &\ar@{-}[ddd]    & \ar@{-}[ddd]&   \ar@{-}[dd]   &  \ar@{-}[ddd] &\\
	&  &   \ar@{-}[rr]  &       & \bullet     \ar@{-}[r] & &  \\
	&   &  \ar@{-}[rrr]&   &  \ar@{-}[d]  &&  \\
	&    &  \ar@{-}[rrr]&   & &&      }
\quad
\xymatrix@R=0.5cm{
	\\
	\\
	\Longrightarrow \\
	\\
	\\}
\quad N=N^*+1,\; \sigma=1
\xymatrix @ = 0.5pc @*[r]{
	\ar@{-}[rrrrr] &     &\ar@{-}[ddd]   &\ar@{-}[ddd]&   \ar@{-}[dd]   &  \ar@{-}[ddd] &\\
	&&   \ar@{-}[rr]  &     & \bullet     \ar@{-}[r] & &  \\
	&     &  \ar@{-}[rrr]&   &  \ar@{-}[d]&&   \\
	&     &  \ar@{-}[rrr]&   & &&     }
$$

\caption{The action of the pfaffian $PfF_{\widehat{2}}$ in the case
$\sigma=0$}\label{fig7}
\end{figure}
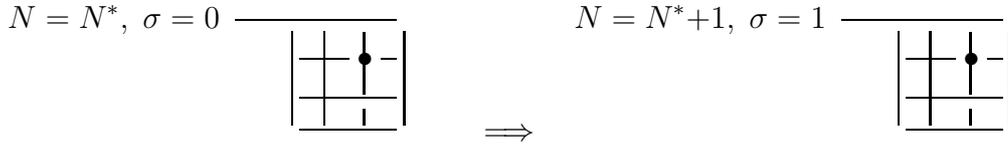

To the formula \ref{pkc2} there corresponds the figure \ref{fig7}

Thus, if one considers the vectors $\xi^{N^*}_t$  for all $t$, then
one gets that the matrix of mapping from the span
$<\xi^{N^*}_t>_{t\in\mathbb{Z}}$ to the span
$<\xi^{N^*+1}_t>_{t\in\mathbb{Z}}$, defined by formula
\ref{pkc1},\ref{pkc2} is the following. If $N=N^*$ is such that
$\sigma=0$, then the matrix is unit. If $N=N^*$ is such that
$\sigma=1$, then the matrix has nonzero $(i,i)$ and $(i+1,i)$.
elements.

At last we can do a step of induction. Consider the cases i)-iii) in
dependence of geometry of the figure which the line $N=N^*$ cuts
from the rectangle on the upper right corner.

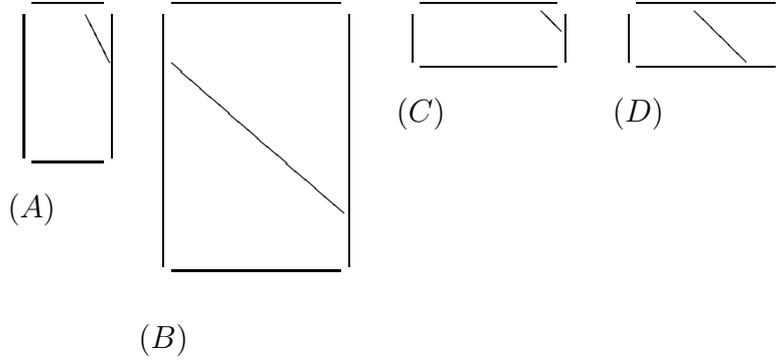
\begin{figure}
\centering
$$
\xymatrix @ = 0.5pc @*[r]{
	\ar@{-}[rr]  \ar@{-}[ddddd]    & \ar@{-}[rdd]   &   \ar@{-}[ddddd]   \\
	\\   &&      \\
	&     &    \\
	&     &  \\
	\ar@{-}[rr] && \\ (A)
} \,\,
\xymatrix@R=0.5cm{
	\ar@{-}[rr]  \ar@{-}[ddddd]   &     &   \ar@{-}[ddddd]   \\
	\ar@{-}[rrddd]\\   &&      \\
	&     &    \\
	&     &  \\
	\ar@{-}[rr] &&\\ (B)
}
\quad \xymatrix @ = 0.5pc @*[r]{
	\ar@{-}[rrrr]  \ar@{-}[dd]   &    && \ar@{-}[rd]  &   \ar@{-}[dd]   \\
	&     & && \\
	\ar@{-}[rrrr] && &&\\(C)
}
\quad \xymatrix @ = 0.5pc @*[r]{
	\ar@{-}[rrrr]  \ar@{-}[dd]   & \ar@{-}[rrdd]   &&   &   \ar@{-}[dd]   \\
	&     & && \\
	\ar@{-}[rrrr] && &&\\(D)
}
$$
\caption{The intersection of the rectangle and the line
$N=const$.}\label{fig8}
\end{figure}

For  the   a narrow rectangle there two different case and  for a
 wide rectangle there  also two case. They are
presented on the figure \ref{fig8}

The considerations in the cases A and C are similar, thus we
actually have three cases.

For a fixed $N=N^*$ we have that for all $N$ either $\sigma=0$ or
$\sigma=1$. Thus in each of three cases above we have to consider
two subcases.

 {\bf Cases A and C.}

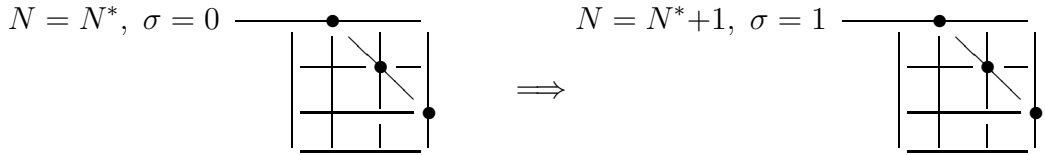
\begin{figure}
\centering
$$N=N^*,\; \sigma=0
\xymatrix @ = 0.5pc @*[r]{
	\ar@{-}[rrrrr] &     &\ar@{-}[ddd]    & \bullet \ar@{-}[ddd] \ar@{-}[rrdd]&   \ar@{-}[dd]   &  \ar@{-}[ddd] &\\
	&  &   \ar@{-}[rr]  &       & \bullet     \ar@{-}[r] & &  \\
	&   &  \ar@{-}[rrr]&   &  \ar@{-}[d]  &\bullet&  \\
	&    &  \ar@{-}[rrr]&   & &&      }
\quad
\xymatrix @ = 0.5pc @*[r]{
	\\
	\\
	\Longrightarrow \\
	\\
	\\}N=N^*+1,\; \sigma=1
\xymatrix @ = 0.5pc @*[r]{
	\ar@{-}[rrrrr] &     &\ar@{-}[ddd]    & \bullet \ar@{-}[ddd] \ar@{-}[rrdd]&   \ar@{-}[dd]   &  \ar@{-}[ddd] &\\
	&  &   \ar@{-}[rr]  &       & \bullet     \ar@{-}[r] & &  \\
	&   &  \ar@{-}[rrr]&   &  \ar@{-}[d]  &\bullet&  \\
	&
	&  \ar@{-}[rrr]&   & &&      }
$$
\caption{The action of the pfaffian $PfF_{\widehat{2}}$ in the case
$\sigma=0$}\label{fig9}
\end{figure}

Firstly consider the case  $\sigma=0$. The pfaffian
$PfF_{\widehat{2}}$ acts as it is shown on the figure \ref{fig9},
i.e. the pfaffian changes $\sigma=0$ to $\sigma=1$.

Obviously in this case the space  with $N=N^*$ is mapped
isomorphically to the space of states with $N=N^*+1$.

Using this isomorphism we assign to the states with $N=N^*+1$ the
fourth quantum number by the following rule. If a state with $N=N^*$
has the fourth quantum number  $k$ then state's image has the fourth
quantum number equal to $k+1$.

By induction one immediately  gets that the one can classify all
states with  $N=N^*+1$ with four quantum numbers.

Consider now the case $\sigma=1$. The pfaffian $PfF_{\widehat{2}}$
acts as it is shown on the figure \ref{fig10}.

\begin{figure}
\centering
$$
\xymatrix @ = 0.5pc @*[r]{
	\ar@{-}[rrrrr] &     & \bullet\ar@{-}[ddd] \ar@{-}[dddrrr]   & \bullet\ar@{-}[ddd]\ar@{-}[ddrr] \ar@{=>}[l]\ar@{=>}[d]&   \ar@{-}[d]   &  \ar@{-}[ddd] &\\
	&&   \ar@{-}[r]  &  \bullet     & \bullet    \ar@{=>}[l]\ar@{=>}[d] \ar@{-}[r] & &  \\
	&     &  \ar@{-}[rrr]&   & \bullet \ar@{-}[d]   &  \bullet \ar@{=>}[l]\ar@{=>}[d] & N=N^*,\; \sigma=1 \\
	&     &  \ar@{-}[rrr]&   &     & \bullet & N=N^*+1,\; \sigma=0 }
$$
\caption{The action of the pfaffian $PfF_{\widehat{2}}$ in the case
$\sigma=1$}\label{fig10}
\end{figure}

Using the formulae of the pfaffian's  and information about the
matrix of this action we see that the pfaffian defines an injective
mapping and the codimension of the image equals $1$.

Choose as  a complementary vector to the image  a vector
corresponding to the lower point on the line  $N=N^*+1$.  To the
corresponding states we assign zero value if the fourth quantum
number. For the states from the image of the pfaffian we use the
following rule. If a state with $N=N^*$ has the fourth quantum
number  $k$, then its image has the fourth quantum number $k+1$.

One gets by induction that   one can classify all states with
$N=N^*+1$ with four quantum numbers.


{\bf Case B.}

\begin{figure}
\centering
$$N=N^*,\; \sigma=0
\xymatrix @ = 0.5pc @*[r]{
	\ar@{-}[rrr]       \ar@{-}[d]    &  \ar@{-}[d] &   \ar@{-}[d]   &  \ar@{-}[d] &\\
	\ar@{-}[rrr] \ar@{-}[ddd]  \bullet  \ar@{-}[rrrddd]  &  \ar@{-}[ddd] &   \ar@{-}[ddd]   &  \ar@{-}[ddd] &\\
	\ar@{-}[rrr]  &     \bullet     &   \ar@{-}[r] & &  \\
	\ar@{-}[rrr]&   &  \bullet &&  \\
	\ar@{-}[rrr]&   & & \bullet&      }
\quad
\xymatrix @ = 0.5pc @*[r]{
	\\
	\\
	\Longrightarrow \\
	\\
	\\} N=N^*+1,\; \sigma=1
\xymatrix @ = 0.5pc @*[r]{
	\ar@{-}[rrr]       \ar@{-}[d]    &  \ar@{-}[d] &   \ar@{-}[d]   &  \ar@{-}[d] &\\
	\ar@{-}[rrr] \ar@{-}[ddd]  \bullet  \ar@{-}[rrrddd]  &  \ar@{-}[ddd] &   \ar@{-}[ddd]   &  \ar@{-}[ddd] &\\
	\ar@{-}[rrr]  &     \bullet     &   \ar@{-}[r] & &  \\
	\ar@{-}[rrr]&   &  \bullet &&  \\
	\ar@{-}[rrr]&   & & \bullet&      }
$$
\caption{The action of the pfaffian $PfF_{\widehat{2}}$ in the case
$\sigma=0$}\label{fig11}
\end{figure}
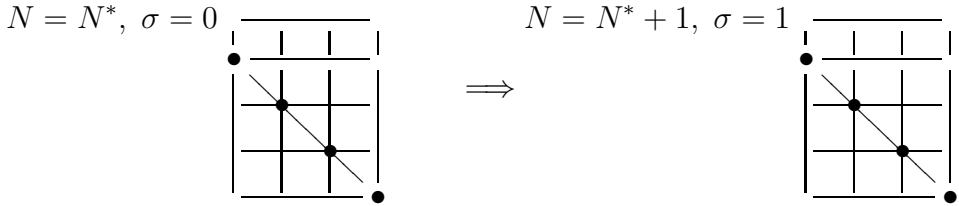

Consider the case when for $N=N^*$ we have $\sigma=0$, then the
pfaffian acts from the space of states with $N=N^*$  to the space of
states with $N=N^*+1$ as it is shown on the figure  \ref{fig11}.

One sees that the pfaffian is an isomorphism of these spaces.

 If a state with $N=N^*$ has the fourth quantum
number  $k$, then we assign to  its image  the fourth quantum number
$k+1$.

One gets by construction that   one can classify all states with
$N=N^*+1$ with four quantum numbers.

If for $N=N^*$ we have $\sigma=1$, then the pfaffian acts from the
space of states with  $N=N^*$ to the spaces of states with $N=N^*+1$
as it is shown on the figure \ref{fig12}.

\begin{figure}
\centering
$$ \xymatrix @ = 0.5pc @*[r]{
	\bullet \ar@{-}[rrr] \ar@{-}[rrrddd]   \ar@{=>}[d]    \ar@{-}[d]    &  \ar@{-}[d] &   \ar@{-}[d]   &  \ar@{-}[d] &\\
	\bullet \ar@{-}[rrr] \ar@{-}[ddd]   \ar@{-}[rrrddd]   & \bullet \ar@{-}[ddd]  \ar@{=>}[l]\ar@{=>}[d] &   \ar@{-}[ddd]   &  \ar@{-}[ddd] &\\
	\ar@{-}[rrr]  & \bullet       & \bullet   \ar@{=>}[l]\ar@{=>}[d]   \ar@{-}[r] & &  \\
	\ar@{-}[rrr]&   &  \bullet  &\bullet  \ar@{=>}[l]\ar@{=>}[d]&  N=N^*,\,\,\, \sigma=1  \\
	\ar@{-}[rrr]&   & &\bullet &     N=N^*+1,\,\,\, \sigma=0  }$$

\caption{The action of the pfaffian $PfF_{\widehat{2}}$ in the case
$\sigma=1$}\label{fig12}
\end{figure}

Again using information about the action of the pfaffian one sees
that the pfaffian is an isomorphism of these spaces.

 If a state with $N=N^*$ has the fourth quantum
number  $k$, then we assign to its image has the fourth quantum
number $k+1$.

One gets by construction that   one can classify all states with
$N=N^*+1$ with four quantum numbers.

Since we consider only the states with $N<0$, then the case of a
wide rectangle is considered completely.

{\bf  Case D }

Let  $N=N^*$ be such that one has $\sigma=1$. Тhen the pfaffian act
 as it shown on the figure \ref{fig13}.

\begin{figure}
\centering
$$N=N^*+1,\, \sigma=0\,\,\,
\xymatrix @ = 0.5pc @*[r]{
	\ar@{-}[rrrr]     \ar@{-}[rrrddd]    \ar@{-}[d]   \bullet&   \ar@{=>}[l]\ar@{=>}[d]  \ar@{-}[d] \bullet    \ar@{-}[rrrddd]    &  \ar@{-}[d] &   \ar@{-}[d]   &  \ar@{-}[d] &\\
	\ar@{-}[rrrr] \ar@{-}[dd]    &   \bullet \ar@{-}[dd] &   \ar@{=>}[l]\ar@{=>}[d]    \bullet \ar@{-}[dd]   &  \ar@{-}[dd] &  \ar@{-}[dd] \\
	\ar@{-}[rrrr] & &     \bullet     &    \ar@{=>}[l]\ar@{=>}[d]   \bullet \ar@{-}[r] & &  \\
	\ar@{-}[rrrr]&  & &  \bullet & \circ & N=N^*,\, \sigma=1  \\
}
$$
\caption{The action of the pfaffian $PfF_{\widehat{2}}$ in the case
$\sigma=1$}\label{fig13}
\end{figure}

One sees that the pfaffian is an injective mapping and the
codimension of the image equals $1$. As  a complimentary vector to
the image we take a vector corresponding to the lower point on the
line $N=N^*+1$.

Assign to this state the zero fourth  number, and to the vectors
from the image of the pfaffian assign the fourth number according to
the following rule.  If a state with $N=N^*$ has the fourth quantum
number  $k$, then its image has the fourth quantum number $k+1$.

One gets by induction that   one can classify all states with
$N=N^*+1$ with four quantum numbers.

\begin{figure}
\centering
$$N=N^*,\; \sigma=0
\xymatrix @ = 0.5pc @*[r]{
	\ar@{-}[rrr] \ar@{-}[ddd]  \bullet  \ar@{-}[rrrddd]  &  \ar@{-}[ddd] &   \ar@{-}[ddd]   &  \ar@{-}[ddd] &\\
	\ar@{-}[rrr]  &     \bullet     &   \ar@{-}[r] & &  \\
	\ar@{-}[rrr]&   &  \bullet &&  \\
	\ar@{-}[rrr]&   & & \bullet&      }
\quad
\xymatrix@R=0.5cm{
	\\
	\\
	\Longrightarrow \\
	\\
	\\} N=N^*+1,\; \sigma=1
\xymatrix @ = 0.5pc @*[r]{
	\ar@{-}[rrr] \ar@{-}[ddd]  \bullet  \ar@{-}[rrrddd]  &  \ar@{-}[ddd] &   \ar@{-}[ddd]   &  \ar@{-}[ddd] &\\
	\ar@{-}[rrr]  &     \bullet     &   \ar@{-}[r] & &  \\
	\ar@{-}[rrr]&   &  \bullet &&  \\
	\ar@{-}[rrr]&   & & \circ&      }
$$
\caption{The action of the pfaffian $PfF_{\widehat{2}}$ in the case
$\sigma=1$}\label{fig14}
\end{figure}
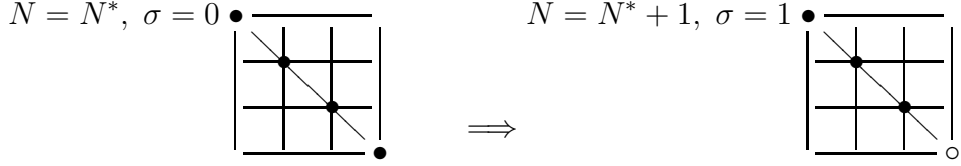

Now let  $N=N^*$ be such that one has $\sigma=0$. Then the pfaffian
acts as it is shown on the figure \ref{fig14}.

The pfaffian  acts as a non-injective mapping, it has a
one-dimensional kernel. This kernel is generated by the state
corresponding to the lower point on the line $N=N^*$. This is just
the vector to which we have assigned a definite quantum number on
the previous step. The image is factor space by the kernel, but
since the vectors in the kernel have definite fourth quantum number
one can  define correctly the fourth quantum number for the vectors
in the factor space. As usual if a state with $N=N^*$ has the fourth
quantum number $k$, then its image has the fourth quantum number
$k+1$.

One gets by construction that   one can classify all states with
$N=N^*+1$ with four quantum numbers.

Thus we have assigned the fourth quantum number for all states with
$N\leq 0$. The induction is completed.

Now let us construct the fourth quantum number for the states with
$N>0$.

There exists a reflection $\omega$ in the Weil group of the root
system $B_2$, that transforms $x$ into $-x$. This element acts  on
the algebra $U(\mathfrak{o}_5)$,  and maps $F_{ij}$ into $F_{ji}$,
one has  $$\omega(PfF_{\widehat{-2}})=-PfF_{\widehat{2}}.$$

Also  $\omega$  acts in every representation $V$ of the algebra
$\mathfrak{o}_5$. Under the action of $\omega$ in  a representation
a weight space
 $V_{\lambda}$ is mapped into  $V_{-\lambda}$.

  The
action of $\omega$ on the algebra and in the representation are
related by the formula
$$\omega(g)\omega(v)=\omega(gv),$$ where $g\in U(\mathfrak{o}_5)$, $v\in
V$.

Now let us construct for the states with $N\geq 0$ the fourth
quantum number by the formula
$$k(v)=k(\omega v).$$

Note that for the states with  $N=0$ one actually obtains two ways
for definition of the fourth quantum number. Fist time we have
assigned the quantum number when we have considered states with
$N\geq 0$ and second time we have assigned the quantum number when
we have considered states with $N\leq 0$. We are going to prove that
they are equivalent.

This follows from the fact the the fourth quantum number together
with  $N$ give an indexation of $\mathfrak{o}_3$-highest vectors
with  $\mathfrak{o}_3$-weight  $T$. If one considers such
$\mathfrak{o}_3$-highest vectors with  $N=0$, then the mapping
$\omega$ acts on the space of such $\mathfrak{o}_3$-highest as
multiplication by scalar. Indeed, since being restricted on this
space $\omega$ turns into analogous $\mathfrak{o}_3$-reflection,
then it preserves $\mathfrak{o}_3$ representations.  Thus each
$\mathfrak{o}_3$-highest vector is mapped into itself multiplied by
a scalar. Since this is true for all $\mathfrak{o}_3$-highest
vectors these scalars are the same for all $\mathfrak{o}_3$-highest
vectors.

That is why two way of construction of the fourth quantum number for
the states with $N=0$ are equivalent.

\endproof

Thus for $N<0$ the pfaffian $PfF_{\widehat{-2}}$ acts as a raising
operator for the quantum number $k$, and for $N>0$ the pfaffian
$PfF_{\widehat{2}}$ acts as a raising operator for the the quantum
number $k$

The theorem  \ref{teork} is proved.

Note that the numbers  $N$ and $T$ give an indexation of base
vectors in the space of $\mathfrak{o}_3$-highest vectors with the
$\mathfrak{o}_3$-weight $T$.

\section{Appendix. Proof of the theorem \ref{maint}}

To prove the theorem we need first to investigate the action of
pfaffians in representations.

\subsection{Action of a pfaffian on a weight vector.}
\label{wei}
 Remind that $e_i$ denotes the standard base vectors $F_{ii}^*$ of  dual space to the Cartan subalgebra.

\begin{prop}
 \label{wes}
Let $V$ be a representation of $\mathfrak{o}_N$.
  Under the action of the pfaffian $PfF_I$  a weight vector  with the weight $\mu$ is mapped to a weight vector with the weight $\mu-\sum_{i\in I}e_i$.
 \end{prop}
\proof
 If  $v$ is a weight vector in a representation  of $\mathfrak{o}_N$ with the weight
  $\mu$, $g_{\alpha}$ is a root element in $\mathfrak{o}_N$ corresponding to the root $\alpha$,
 then $g_{\alpha}v_{\mu}$ is a weight vector of to the weight  $\alpha+\mu$.

Consider the vector $PfF_Iv$.
  By definition one has
  $$PfF_I=\frac{1}{\frac{k}{2}!2^{\frac{k}{2}}}\sum_{\sigma\in
S^k}(-1)^{\sigma}F_{-\sigma(i_1)\sigma(i_2)}...F_{-\sigma(i_{k-1})\sigma(i_k)}.$$

To prove the proposition it suffices to show that every summand
changes the weight by subtracting of the same expression
$-\sum_{i\in I}e_i$. Using the correspondence between roots and
elements $F_{ij}$ from the sec. \ref{or} one gets the following.
When one acts by
$F_{-\sigma(i_1)\sigma(i_2)}...F_{-\sigma(i_{k-1})\sigma(i_k)}$ on
$v$ then to the weight the vector
$$e_{-\sigma(i_1)}-e_{\sigma(i_2)}-...+e_{-\sigma(i_{k-1})}-e_{\sigma(i_k)}=-\sum_{i\in
I}e_i$$ is added. This proves the proposition.

 \endproof

Consider the most interesting case
$\mathfrak{o}_{N}=\mathfrak{o}_{2n+1}$ and $|I|=2n$.

\begin{cor}
\label{wes1} Let $\mathfrak{o}_{N}=\mathfrak{o}_{2n+1}$.

The action of $PfF_{\widehat{-n}}$ adds the vector
 $-\sum_{i\in I}e_i=-e_{n}$ to the weight.

The action of $PfF_{\widehat{n}}$ adds the vector
 $-\sum_{i\in I}e_i=-e_{-n}=e_{n}$ to the weight.
\end{cor}

\subsection{Commutators of pfaffians and $F_{ij}$.}

\label{com}

\begin{lem}
\label{l2}
  Let $I=\{i_1,...,i_k\}$, where $k$ is even. Then the commutator $[PfF_I,F_{j_1-j_2}]$
   is calculated according to the following rule.

\begin{enumerate}

\item If $j_1,j_2\notin I$, then  $[PfF_I,F_{j_1-j_2}]=0$.

\item If $j_1\in I,j_2\notin I$, then $[PfF_I,F_{j_1-j_2}]=PfF_{I_{j_1\rightarrow -j_2}}$.

\item If $j_1\notin I,j_2\in I$, then  $[PfF_I,F_{j_1-j_2}]=-PfF_{I_{j_2\rightarrow -j_1}}$.

\item If $j_1\in I,j_2\in I$, then
  $[PfF_I,F_{j_1-j_2}]=PfF_{I_{j_1\rightarrow -j_2}}-PfF_{I_{j_2\rightarrow -j_1}}$
\end{enumerate}

\end{lem}

\proof

The quadratic form is $G=(\delta_{i,-j})$. Hence one can identify
$E_{ij}$ with $e_{i}\otimes e_{-j}$. Then $F_{ij}$ is identified
with $e_{i}\wedge e_{-j}$. Remind that
 $$Pf
F_I=\frac{1}{2^{\frac{k}{2}}\frac{k}{2}!}\sum_{ \sigma\in
S_k}(-1)^{\sigma}F_{-\sigma(i_1)\sigma(i_2)}...F_{-\sigma(i_{k-1})\sigma(i_{k})}.$$

Thus  $PfF_I$ with indexing set $I=\{i_1,...,i_k\}$ is  identified
with the polyvector $e_{-i_1}\wedge...\wedge e_{-i_k}$.

This identification is compatible with the action of
$\mathfrak{o}_N$. Thus

\begin{center}
$[PfF_I,F_{j_1-j_2}]=-[F_{j_1-j_2},PfF_I]=F_{j_1-j_2}e_{i_1}\wedge
e_{i_2}...\wedge e_{i_k}+e_{i_1}\wedge F_{j_1-j_2}e_{i_2}... \wedge
e_{i_k}+e_{i_1}\wedge e_{i_2}... \wedge
F_{j_1-j_2}e_{i_k}.$\end{center}

One has $F_{j_1-j_2}e_{-i_1}=e_{j_2}$ if $j_2=i_1$ and
$F_{j_1-j_2}e_{-i_1}=-e_{j_2}$ if $j_1=i_1$.

Suggest that $\{j_1,j_2\}\cap I=\emptyset$.   Then
$[PfF_I,F_{j_1-j_2}]=0$.

Suggest that $j_1\in I,j_2\notin I$. Then  $j_1=i_t$ and the only
nonzero summand is that containing $F_{j_1j_2}e_{-i_t}$.  Thus we
have $[PfF_I,F_{j_1-j_2}]=-PF_{I_{j_1\mapsto -j_2}}$. Here
$I_{j_1\mapsto -j_2}$ is obtained from $I$ by replacing the index
$j_1$ to $j_2$.

The case $j_1\notin I,j_2\in I$ is considered in the same manner.

Suggest that $\{j_1,j_2\}\subset I$. That is $j_1=e_{i_t}$,
$j_2=e_{i_s}$. Then in the sum there are two nonzero summands one
contains  $F_{j_1j_2}e_{i_t}$, the other contains
$F_{j_1j_2}e_{i_s}$.  Each of them is a wedge product with a new
indexing set. So one gets that
$[PfF_I,F_{j_1j_2}]=PfF_{I_{j_1\rightarrow
-j_2}}-PfF_{I_{j_2\rightarrow -j_1}}$.

\endproof

\begin{cor}
\label{cor3} In the case $\mathfrak{o}_{2n+1}$  the pfaffians
$PfF_{\widehat{n}}$, $PfF_{\widehat{-n}}$ commute with elements
$F_{ij}$, $-n<i,j<n$, that span the subalgebra $\mathfrak{o}_{2n-1}$

\end{cor}

\subsection{Some formulas involving pfaffians.}
\label{formu}

In this subsection some  summation formulas  are proved.

\begin{lem}
\label{minors1}

$PfF_I=\frac{(\frac{p}{2})!(\frac{q}{2})!}{(\frac{k}{2})!}\sum_{I=I'\sqcup
I'',|I'|=p,|I''|=q}(-1)^{(I'I'')}PfF_{I'}PfF_{I''}$.

Here $(-1)^{(I'I'')}$ is a sign of a permutation of the set
$I=\{i_1,...,i_k\}$ that places first the subset $I'\subset I$ and
then the subset $I''\subset I$.

The numbers $p$, $q$ are even fixed numbers,  they satisfy
$p+q=k=|I|$.
\end{lem}
\proof

By definition one has

$$PfF_I=\frac{1}{2^{\frac{k}{2}}(\frac{k}{2})!}\sum_{\sigma\in
S_k}(-1)^{\sigma}F_{-\sigma(i_1)\sigma(i_2)}...F_{-\sigma(i_{k-1}),\sigma(i_k)}.$$

The summand $
(-1)^{\sigma}F_{-\sigma(i_1)\sigma(i_2)}...F_{-\sigma(i_{k-1}),\sigma(i_k)}$
can be written as

$$(-1)^{(I'I'')}(-1)^{\sigma'}F_{-\sigma'(i'_1)\sigma'(i'_2)}...F_{-\sigma'(i'_{p-1}),\sigma'(i'_p)}
(-1)^{\sigma''}F_{-\sigma''(i''_1)\sigma''(i''_2)}...F_{-\sigma''(i''_{q-1}),\sigma''(i''_q)}.$$

Here $I'=\{i'_1,...,i'_p\}$ is the set of indices
$\{\sigma(i_1),...,\sigma(i_p)\}$ placed in a natural order,
$I''=\{i''_1,...,i''_q\}$ is set   of indices
$\{\sigma(i_{p+1}),...,\sigma(i_k)\}$ placed in a natural order,
$\sigma'$ is a permutation $\{\sigma(i_1),...,\sigma(i_p)\}$  of the
set $I'$ and  $\sigma''$ is a permutation of the set $I''$ defined
in a similar way. Note that
$(-1)^{(I'I'')}(-1)^{\sigma'}(-1)^{\sigma''}=(-1)^{\sigma}$.

The mapping $\sigma\mapsto I',I'',\sigma',\sigma''$ is bijective.

Thus the pfaffian can be written as

$\frac{(\frac{p}{2})!(\frac{q}{2})!}{(\frac{k}{2})!}\sum_{I=I'\sqcup
I'',|I'|=p,|I''|=q}(-1)^{(I'I'')}\frac{1}{2^{\frac{k}{2}}(\frac{p}{2})!(\frac{q}{2})!}\sum_{\sigma'}(-1)^{\sigma'}(-1)^{\sigma''}
F_{-\sigma'(i'_1)\sigma'(i'_2)}...\\...F_{-\sigma'(i'_{p-1}),\sigma'(i'_p)}F_{-\sigma''(i''_1)\sigma'(i''_2)}...F_{-\sigma''(i''_{q-1}),\sigma''(i''_q)}=\\
=\frac{(\frac{p}{2})!(\frac{q}{2})!}{(\frac{k}{2})!}\sum_{I=I'\sqcup
I'',|I'|=p,|I''|=q}(-1)^{(I'I'')}PfF_{I'}PfF_{I''}$
\endproof

\begin{cor}
\label{corl2} If $|I|=k$, then

$PfF_I=\frac{1}{\frac{k}{2}+1}\sum_{I=I'\sqcup
I''}(-1)^{(I'I'')}\frac{(\frac{|I'|}{2})!(\frac{|I''|}{2})!}{(\frac{k}{2})!}PfF_{I'}PfF_{I''}$.

\end{cor}

\begin{lem}
\label{minorn}

 Let $-n\in I$. Then
$PfF_I=\\=\sum_{i\in I\setminus\{-n\}}\sum_{I\setminus
\{-n,i\}=I'\sqcup
I''}\frac{(\frac{|I'|}{2})!(\frac{|I''|}{2})!}{(\frac{k}{2})!}(-1)^{(I'-niI'')}PfF_{I'}F_{ni}PfF_{I''}$.

Here $(-1)^{(I'-niI'')}$ is a sign of the permutation $(I',-n,i,I'')$
of the set $I$.
\end{lem}
\proof

By definition one has

$$PfF_I=\frac{1}{2^{\frac{k}{2}}(\frac{k}{2})!}\sum_{\sigma\in
S_k}(-1)^{\sigma}F_{-\sigma(i_1)\sigma(i_2)}...F_{-\sigma(i_{k-1}),\sigma(i_k)}.$$

Since $F_{ij}=-F_{-j-i}$ the summation can be taken only over such
permutation  such that $\sigma(i_{2t-1})<\sigma(i_{2t})$. But if the
summation is done in such a way the multiple
$\frac{1}{2^{\frac{k}{2}}}$ must be omitted.

Fix a such a permutation $\sigma$ and find a place such that
$(\sigma(i_{2t-1}),\sigma(i_{2t}))=(-n,i)$. The summand
$(-1)^{\sigma}F_{-\sigma(i_1)\sigma(i_2)}...F_{-\sigma(i_{k-1}),\sigma(i_k)}$
can be written as

$(-1)^{(I'inI'')}(-1)^{\sigma'}F_{-\sigma'(i'_1)\sigma'(i'_2)}...F_{-\sigma'(i'_{p-1}),\sigma'(i_p)}F_{-ni}(-1)^{\sigma''}F_{-\sigma''(i''_1)\sigma''(i''_2)}...
F_{-\sigma''(i''_{q-1}),\sigma''(i''_q)}$.

Here $I'=\{i'_1,...,i'_p\}$ is the set of indices
$\{\sigma(i_1),...,\sigma(i_{2t-2})\}$ placed in the natural order,
$I''=\{i''_1,...,i''_q\}$ is set   of indices
$\{\sigma(i_{2t+1}),...,\sigma(i_k)\}$ placed in the natural order,
$\sigma'$ is a permutation $\{\sigma(i_1),...,\sigma(i_{2t-2})\}$ of
the set $I'$ and $\sigma''$ is a permutation of the set $I''$
defined in a similar way. Note that
$(-1)^{(I'-niI'')}(-1)^{\sigma'}(-1)^{\sigma''}=(-1)^{\sigma}$. The
permutation $\sigma'$  satisfies the condition
$\sigma'(i'_{2t-1})<\sigma'(i'_{2t})$ as well as  the permutation
$\sigma''$.

The mapping $\sigma\mapsto I',I'',\sigma',\sigma''$ is bijective
(since $\sigma(i_{2t-1})<\sigma(i_{2t})$).

Thus the pfaffian can be written as

$\sum_{i\in I\setminus\{n\}}\sum_{I\setminus \{i,n\}=I'\sqcup I''}
\frac{(\frac{|I'|}{2})!(\frac{|I''|}{2})!}{(\frac{k}{2})!}\frac{1}{(\frac{|I'|}{2})!(\frac{|I''|}{2})!}(-1)^{(I'-niI'')}
(-1)^{\sigma'}F_{-\sigma'(i'_1)\sigma'(i'_2)}...\\...F_{-\sigma'(i'_{p-1}),\sigma'(i_p)}F_{ni}(-1)^{\sigma''}F_{-\sigma''(i''_1)\sigma''(i''_2)}...
F_{-\sigma''(i''_{q-1}),\sigma''(i''_q)}=\\=\sum_{i\in
I\setminus\{n\}}\sum_{I\setminus \{i,n\}=I'\sqcup
I''}\frac{(\frac{|I'|}{2})!(\frac{|I''|}{2})!}{(\frac{k}{2})!}(-1)^{(I'-niI'')}PfF_{I'}F_{ni}PfF_{I''}$

\endproof

\subsection{The Mickelsson-Zhelobenko algebra and
pfaffians.}\label{mj}

In the construction of the Gelfand-Tsetlin-Molev base a key role is
played by the  Mickelsson-Zhelobenko algebra.

 There exist a
projection from the algebra $U(\mathfrak{o}_N)$  to the
Mickelsson-Zhelobenko algebra
$Z(\mathfrak{o}_{2n+1},\mathfrak{o}_{2n-1})$. To prove the theorem
\ref{maint} we must calculate the image of pfaffians under this
projection. The present section is devoted to the calculation of
this image.

\subsubsection{The Mickelsson-Zhelobenko algbera}\label{mj1}

Remind the idea of the construction of the Gelfand-Tsetlin-Molev
base of  a representation $V$ of the algebra $\mathfrak{o}_{2n+1}$.
An irreducible  representation  $V$  of the algebra
$\mathfrak{o}_{2n+1}$ becomes reducible as a representation
$\mathfrak{o}_{2n-1}$. According to the idea of Gelfand and Tsetlin
to construct a base in $V$, it is necessary to know all possible
highest weights $\mu$ of  $\mathfrak{o}_{2n-1}$-irreps, into which
$V$ splits. Also it is necessary to be able to construct a base in
the multiplicity space, that is in the space of
$\mathfrak{o}_{2n-1}$-highest vectors with the fixed
$\mathfrak{o}_{2n-1}$-weight $\mu$.

The first problem is solved quite easily, for the solution of the
second one the Mickelsson-Zhelobenlo algebra is used.

\begin{defn}Denote as $V_{\mu}^{+}$ the space of $\mathfrak{o}_{2n-1}$-highest vectors with the
$\mathfrak{o}_{2n-1}$-weight  $\mu$ in the
$\mathfrak{o}_{2n+1}$-representation $V$.

\end{defn}

To construct a base in $V_{\mu}^{+}$ Molev has used the
Mickelsson-Zhelobenko algebra, acting on the space
$\oplus_{\mu}V_{\mu}^+$.

Let us give the definition of the Mickelsson-Zhelobenko algebra, see
 \cite{JL},\cite{Molev1}. Al facts and notations are borrowed from   and \cite{Molev}.

Let $\mathfrak{g}$ be a Lie algebra and let $\mathfrak{k}$ be its
reductive subalgebra. The main example is
$\mathfrak{g}=\mathfrak{o}_{2n+1}$ and
$\mathfrak{k}=\mathfrak{o}_{2n-1}$. Let
$\mathfrak{k}=\mathfrak{k}^-+\mathfrak{h}+\mathfrak{k}^+$ be a
triangular decomposition. Let $R(\mathfrak{h})$ be a field of
fractions of the algebra $U(\mathfrak{h})$. Denote
$$U'(\mathfrak{g})=U(\mathfrak{g})\otimes_{U(\mathfrak{h})}
R(\mathfrak{h}).$$

Let $$J'=U'(\mathfrak{g})\mathfrak{k}^+$$ be the left ideal in
$U'(\mathfrak{g})$, generated by $\mathfrak{k}^+$. Put
$$M(\mathfrak{g},\mathfrak{k})=U'(\mathfrak{g})/J'.$$

For every positive root $\alpha$ of the algebra $\mathfrak{k}$
define a series

$$p_{\alpha}=1+\sum_{k=1}^{\infty}e_{-\alpha}^ke_{\alpha}^k\frac{(-1)^k}{k!(h_{\alpha}+\rho(h_{\alpha})+1)...(h_{\alpha}+\rho(h_{\alpha})+k)},$$

where $e_{\alpha}$ is a root vector  $\mathfrak{k}$, corresponding
to $\alpha$, $h_{\alpha}$ is a corresponding Cartan element, $\rho$
is a half-sum of positive roots of $\mathfrak{k}$.

An order is normal if the following holds. Let a root
  be
   a sum of two roots, then it lies between them. Chose a normal ordering $\alpha_1<...<\alpha_m$  of  positive roots
of
  $\mathfrak{k}$.

Put $$p=p_{\alpha_1}...p_{\alpha_m}.$$ This element is called the
extremal projector. It can be proved that nevertheless  $p$ is an
infinite series it's action on $M(\mathfrak{g},\mathfrak{k})$  by
left multiplication is well defined \cite{JL}.

The following equalities hold: $e_{\alpha}p=pe_{-\alpha}=0$, where
$\alpha$  is a positive root of $\mathfrak{k}$.

Put $$Z(\mathfrak{g},\mathfrak{k})=pM(\mathfrak{g},\mathfrak{k}).$$
This is the Mickelsson-Zhelobenko algebra. The multiplication in
$Z(\mathfrak{g},\mathfrak{k})$ is defined using the isomorphism
$Z(\mathfrak{g},\mathfrak{k})=Norm J'/J'$, where $NormJ'=\{u\in
U'(\mathfrak{g}): J'u\subset J'\}$. Thus
$Z(\mathfrak{g},\mathfrak{k})$ is an associative algebra and a
bimodule over $R(\mathfrak{h})$ \cite{JL}.

There exists a natural projection from $U(\mathfrak{o}_N)$  to
$Z(\mathfrak{o}_{2n+1},\mathfrak{o}_{2n-1})$, which sends $x\in
U(\mathfrak{o}_N)$ to $p\,x\,mod J'$.

Choose linear independent elements $v_1,...,v_n\in \mathfrak{g}$,
 such that  $<v_1,...,v_n>\oplus\mathfrak{ k}=\mathfrak{g}$ as linear spaces over  $\mathbb{C}$. Put $z_i=pv_i mod J'$.
It can be proved that monomials
$\check{z}_1^{m_1}...\check{z}_n^{m_n}$, $m_i\in\mathbb{Z}^+$, form
a bases of $Z(\mathfrak{g},\mathfrak{k})$ over $R(\mathfrak{h})$.

In the case $Z(\mathfrak{o}_{2n+1},\mathfrak{o}_{2n-1})$  put
$$\check{z}_{i\pm n}=pF_{i,\pm n}modJ', \,\,\,\,i=-n,...,n.$$ Notations are taken from \cite{Molev}. There exists an obvious
symmetry $\check{z}_{ij}=\check{z}_{-j-i}$. From previous
considerations  it follows that
$Z(\mathfrak{o}_{2n+1},\mathfrak{o}_{2n-1})$  is generated by
elements  $\check{z}_{ia}$, $i=0,...,n$, $a=\pm n$ or
$\check{z}_{ai}$, $i=0,...,n$, $a=\pm n$.

Sometimes it is more useful to use the generators

$$z_{i\pm n}=\check{z}_{i\pm n}(f_i-f_{i-1})...(f_i-f_{-n+1}),$$

where $$f_i=F_{ii}+ \rho_i, \text{ for $i>0$
},\,\,\,f_0=-\frac{1}{2},\,\,\, f_{-i}=-f_i,$$ and
$$\rho_i=i-\frac{1}{2} \text{ for $i>0$ and } \rho_{-i}=-\rho_{i}.$$

 In particular

 $$z_{0n}=\check{z}_{0n}\prod_{i=1}^{n-1}(F_{ii}+i-\frac{1}{2}).$$

The Mickelsson-Zhelobenko algebra
$Z(\mathfrak{o}_{2n+1},\mathfrak{o}_{2n-1})$  acts on the space
$\oplus_{\mu}V_{\mu}^+$ (see \cite{Molev}). A weight $\mu$ changes
under this action according to the following rule. Let $a$ be $\pm
n$ and
$\mu+\delta_i=(\mu_1,...,\mu_{i-1},\mu_i+1,\mu_{i+1},...,\mu_{n-1})$.
Then for $i=1,...,n-1$ the following holds

$$z_{ia}:V_{\mu}^+\rightarrow V_{\mu+\delta_i}^{+},\,\,\,\,z_{ai}:V_{\mu}^+\rightarrow V_{\mu-\delta_i}^{+}$$

Elements $z_{0a}$ do not change a $\mathfrak{o}_{2n-1}$-weight, that
is they map $V_{\mu}^+$ into itself.




\subsubsection{Images of pfaffians in the Mickelsson-Zhelobenko
alebra } \label{mjp1}


\begin{defn}
 A product  of root and Cartan elements in the universal enveloping
  algebra is called normally ordered if in it at first (from the left side)  the negative root elements
   occur, then Cartan elements occur and at the end  positive root elements occur.
\end{defn}

  Every product of  root and Cartan elements  equals to a sum of normally ordered products.


\begin{prop} Let  $I\subset\{-n+1,...,n-1\}$ be a subset which is not symmetric with respect to zero. Then $pPfF_I=0$ in
$Z(\mathfrak{o}_{2n+1},\mathfrak{o}_{2n-1})$ or in
$Z(\mathfrak{o}_{2n},\mathfrak{o}_{2n-2})$.
\end{prop}
 \proof
According to the definition a pfaffian a sum over permutations.  The
summands  are products of root vectors and Cartan elements  of
$\mathfrak{o}_N$

 The sum of root corresponding to element of each product equals
 $-\sum_{i\in I}e_i$. Since the set $I$ is nonsymmetric one has $-\sum_{i\in I}e_i\neq 0$.

Impose a normal ordering in every summand. When one does the normal
ordering new summands appear. But from the equality
$[e_{\alpha},e_{\beta}]=N_{\alpha,\beta}e_{\alpha+\beta}$ it follows
that the sum of roots corresponding to the elements of these new
products is again  $-\sum_{i\in I}e_i$.

Since $-\sum_{i\in I}e_i\neq 0$  in every normally ordered summand in the pfaffian there is a root element.
 These elements  either are zero modulo $J'$,  if there is  a positive root element, or  vanish  after multiplication by $p$,
  if there is a negative root element.

\endproof

Let us give a formula for the image of a pfaffian whose indexing set
$I$ is symmetric and is contained in $\{-n+1,...,n-1\}$. In this
case the calculation of the image in the Mickelsson-Zhelobenko
algebra is equivalent to the calcualtion of the image of the
pfaffian under the Harish-Candra homomorphism. This calculation was
done in \cite{ce2} (proposition 7.1), the result is the following.

\begin{prop}\cite{ce2}
\label{tps}
$PfF_I=\frac{1}{(\frac{k}{2})!}D_{\frac{k}{2}}(F_{i_1i_1},...,F_{i_{\frac{k}{2}\frac{k}{2}}}),$

where $D_r(h_1,...,h_r)=\prod_{i=1}^r(h_i-\frac{r}{2}+i)$

\end{prop}


Now let us find an image in
$Z(\mathfrak{o}_{2n+1},\mathfrak{o}_{2n-1})$ of the pfaffian
$PfF_{\widehat{n}}$.

To formulate the next theorem define a polynomial $C_n$.

\begin{defn}
\label{ch} Let
$C_{n-1}(h_1,...,h_{n-1})=(-1)^{n-1}D_{n-1}(h_1,...,h_{n-1})-\\-4\sum_{i=1}^{n-1}(-1)^{t+1}D_{n-2}(h_1,...,\widehat{h_i},...,h_{n-1})$

\end{defn}

\begin{thm}
\label{tm} The image of $PfF_{\widehat{n}}$ in
$Z(\mathfrak{o}_{2n+1},\mathfrak{o}_{2n-1})$
 equals $\check{z}_{n0}C_{n-1}(F_{11},...,F_{(n-1)(n-1)})$.

\end{thm}
\proof

Take a set of indices of type
$I=\{-n,-i_{\frac{k-2}{2}},...,-i_{1},0,i_1,...,i_{\frac{k-2}{2}}\}$



By Lemma \ref{minorn} the following equality takes place

$$PfF_{I}=\sum_{i\in I\setminus \{-n\}}\sum_{I'\sqcup
I''=I\setminus\{i,-n\}}\frac{(\frac{|I'|}{2})!(\frac{|I''|}{2})!}{(\frac{k}{2})!}(-1)^{(I'-niI'')}PfF_{I'}F_{ni}PfF_{I''}.$$

To find  the image in the Mickelson-Zhelobenko algebra of the sum
$\sum_{I'\sqcup I''=I\setminus\{i,-n\}}$ divide the summands into
three groups: $1)$ those for which $i=0$; $2)$ those for which
$i<0$; $3)$ those for which $i>0$.

Let us find the  image of summands for which $i=0$. In this case
$PfF_{I'}$ and $PfF_{I'}$ commute with $F_{n0}$. Note that
$(-1)^{(I'-n0I'')}=(-1)^{(I'I'')}(-1)^{\frac{k}{2}-1}$ (to prove
this firstly move $-n,0$  to two first places and then move $0$ to
the right place, then the signs $(-1)^{|I'|}$, $(-1)^{|I'|-1}$,
$(-1)^{\frac{k}{2}}$ appear).

Using the corollary \ref{corl2} one gets that the image sum these
summands equals
\begin{center}
$(-1)^{\frac{k}{2}-1}\sum_{I'\sqcup
I''=I\setminus\{-n,0\}}\frac{(\frac{|I'|}{2})!(\frac{|I''|}{2})!}{(\frac{k}{2})!}(-1)^{(I'I'')}PfF_{I'}PfF_{I''}F_{n0}=
\frac{1}{\frac{k}{2}}(-1)^{\frac{k}{2}-1}PfF_{I\setminus\{-n0\}}F_{n0}=(-1)^{\frac{k}{2}-1}PfF_{I\setminus\{-n,0\}}F_{n0}.$
\end{center}

Since the sets of indices $\pm (I\setminus\{-n0\})$ and $\{-n,0\}$
do not intersect, one can apply the projector $p$ and equivalence
$modJ'$ to each multiple.

Thus  the image of these summands is
$$\check{z}_{n0}(pPf_{I\setminus\{-n0\}}mod J').$$

Find the image of summands
$$\frac{(\frac{|I'|}{2})!(\frac{|I''|}{2})!}{(\frac{k}{2})!}(-1)^{(I'-niI'')}PfF_{I'}F_{ni}PfF_{I''},$$
for which $i\neq 0$. Let  $i>0$. Then change $F_{ni}$ and
$PfF_{I''}$. One obtains an expression
$$\frac{(\frac{|I'|}{2})!(\frac{|I''|}{2})!}{(\frac{k}{2})!}(-1)^{(I'-niI'')}(PfF_{I'}PfF_{I''}F_{ni}-PfF_{I'}[PfF_{I''},F_{ni}]).$$
 Now let $i<0$.
Change $F_{ni}$ and $PfF_{I'}$,
 one gets

 $$\frac{(\frac{|I'|}{2})!(\frac{|I''|}{2})!}{(\frac{k}{2})!}(-1)^{(I'-niI'')}(F_{ni}PfF_{I'}PfF_{I''}-[PfF_{I'},F_{ni}]PfF_{I''}).$$

Consider the case  $i>0$. In the last expression the first summand
has a zero image in the Mickelsson-Zhelobenko algebra by the
following reason.
 The sum of roots corresponding to the elements $F_{ij}$ that
 participate in the expression for
$PfF_{I'}F_{ni}PfF_{I''}$ equals to $e_{n}$. The element $F_{ni}$
corresponds to the root $e_{n}-e_{i}$. Thus the sum of roots
corresponding to the elements  $PfF_{I'}PfF_{I''}$ equals $-e_i$.
 Express $PfF_{I'}PfF_{I''}$ as a sum of normally ordered products. Since $i>0$ than in every obtained
   normally product there is a negative root element of the algebra  $\mathfrak{o}_{2n-1}$. Thus after applying the extremal projector $p$ the expression
$PfF_{I'}PfF_{I''}$ vanishes.

In the case $i<0$ it is similarly proved that the first summand has
a zero image in the Mickelsson-Zhelobenko algebra.

Now consider the second summand $$-PfF_{I'}[PfF_{I''},F_{ni}]$$ in
the case $i>0$ or $$-[PfF_{I'},F_{ni}]PfF_{I''}$$ in the case $i<0$.

In the  case $i>0$ if $-i\notin I''$ the expression is zero and
otherwise  it equals to $$-PfF_{I'}PfF_{I''}\mid_{-i\mapsto -n}.$$
 In the case $i<0$    if $-i\notin I'$ the expression is zero and  otherwise it
equals  to $$-PfF_{I'}\mid_{-i\mapsto -n}PfF_{I''}.$$

Thus the image of summands for which $i\neq 0$ equals to the image
of the expression
\begin{center}
$ -\sum_{i\in I\setminus \{-n\},i>0}\sum_{I'\sqcup
I''=I\setminus\{-n,i\},-i\in
I''}\frac{(\frac{|I'|}{2})!(\frac{|I''|}{2})!}{(\frac{k}{2})!}(-1)^{(I'-niI'')}PfF_{I'}PfF_{I''}\mid_{-i\mapsto
-n}- \sum_{i\in I\setminus \{-n\},i<0}\sum_{I'\sqcup
I''=I\setminus\{-n,i\},-i\in
I'}\frac{(\frac{|I'|}{2})!(\frac{|I''|}{2})!}{(\frac{k}{2})!}(-1)^{(I'-niI'')}PfF_{I'}\mid_{-i\mapsto
-n}PfF_{I''}$\end{center}

Let us prove a proposition.

\begin{prop}
The expression above equals
 \begin{equation}\label{gy}-2\sum_{t=-\frac{k}{2},\neq
0}^{\frac{k}{2}}(-1)^{\frac{k}{2}-t-1}\sum_{J'\sqcup
J''=I\setminus\{\pm
i\}}\frac{(\frac{|J'|}{2})!(\frac{|J''|}{2})!}{(\frac{k}{2})!}(-1)^{(J'J'')}PfF_{J'}PfF_{J''}.\end{equation}

\end{prop}

\proof To prove this let us  firstly calculate the sign
$(-1)^{(I'-niI'')}$. The sign $(-1)^{(I'-niI'')}$  differs from the
sign $(-1)^{(I'I'')}$ by the sign of the permutation which moves
$-n,i$ to their right places. This permutation can be done as
follows: first of all move $-n,i$ to two last places,  then move $i$
to it's right place. If $i=i_t$, then
$$(-1)^{(I'-niI'')}=(-1)^{(I'I'')}(-1)^{(|I'|+|I'|+\frac{k-2}{2}-t)}=(-1)^{\frac{k}{2}-t-1}(-1)^{(I'I'')}.$$

Secondly compare $PfF_{I'}\mid_{-i\mapsto -n}$,
$PfF_{I''}\mid_{-i\mapsto -n}$ and
$PfF_{(I'\setminus\{-i\})\cup\{-n\}}$,
$PfF_{(I''\setminus\{-i\})\cup\{-n\}}$ respectively. Here it is
assumed that $i\in I'$ and $i\in I''$.
 In all these expressions at first, the index  $-i$  is changed to
$n$,  but then in the last two expressions the new set of indices is naturally ordered.
  Thus $PfF_{I'}\mid_{-i\mapsto -n}$ and  $PfF_{(I'\setminus\{-i\})\cup\{-n\}}$, $PfF_{I''}\mid_{-i\mapsto -n}$ and  $PfF_{(I''\setminus\{-i\})\cup\{-n\}}$, differ by the sign of this ordering.

  For summands in the sum $\sum_{i\in I\setminus \{-n\},i<0}\sum_{I'\sqcup
I''=I\setminus\{-ni\},-i\in I''}$ denote
  $$J':=(I'\setminus\{-i\})\cup\{-n\}, \,\,\,J'':=I''.$$ One obtains
  $$(-1)^{(I'I'')}PfF_{I'}\mid_{-i\mapsto
  -n}PfF_{I''}=(-1)^{(J'J'')}PfF_{J'}PfF_{J''}.$$ The sign that appears after the ordering is contained in
  $(-1)^{(J'J'')}$.

  Analogously for the summands in the sum $\sum_{i\in I\setminus \{-n\},i<0}\sum_{I'\sqcup
I''=I\setminus\{-ni\},-i\in I''}$, denote
$$J':=I',\,\,\,J'':=(I''\setminus\{-i\})\cup\{-n\}.$$ One obtains that
$$(-1)^{(I'I'')}PfF_{I'}\mid_{-i\mapsto
  -n}PfF_{I''}=(-1)^{(J'J'')}PfF_{J'}PfF_{J''}.$$

In both cases one has $J'\sqcup J''=I\setminus\{\pm i \}$. Also $|J'|=|I'|$, $|I''|=|J''|$.

Note that a pair of sets $J',J''$ occurs twice.  First  as
$(I'\setminus\{-i\})\cup\{-n\}$, $I''$, second as $I'$,
$(I''\setminus\{-i\})\cup\{-n\}$.

Thus one obtains that the considered sum of images of summands for
which $i\neq 0$ is given by the expression \ref{gy}

 $$2\sum_{t=-\frac{k}{2},\neq
0}^{\frac{k}{2}}(-1)^{\frac{k}{2}-t-1}\sum_{J'\sqcup
J''=I\setminus\{\pm
i_t\}}\frac{(\frac{|J'|}{2})!(\frac{|J''|}{2})!}{(\frac{k}{2})!}(-1)^{(J'J'')}PfF_{J'}PfF_{J''}$$

The proposition is proved.
\endproof

By Corollary \ref{corl2} this expression \ref{gy} equals
$$2\sum_{t=-\frac{k}{2},\neq
0}^{\frac{k}{2}}(-1)^{\frac{k}{2}-t-1}PfF_{I\setminus\{\pm i_t\}}
=4\sum_{t=1}^{\frac{k}{2}}(-1)^{\frac{k}{2}-t-1}PfF_{I\setminus\{\pm
i_t\}}. $$

Finally in $Z(\mathfrak{o}_{2n+1},\mathfrak{o}_{2n-1})$ one has
\begin{equation} \label{eee}
PfF_{I}=(-1)^{\frac{k}{2}-1}\check{z}_{n0}(pPfF_{I\setminus
\{-n,0\}}mod
J')-4\sum_{t=1}^{\frac{k}{2}}(-1)^{\frac{k}{2}-t-1}PfF_{I\setminus\{\pm
i_t\}}\end{equation}

Note that $PfF_{I\setminus\{-n,\pm i_t\}}$  is a pfaffian
$PfF_{I^{t}}$  for a new indexing set $I^t=I\setminus\{\pm i_t\}$.

This set is of the same type as $I$.  Apply to   each pfaffian
$PfF_{I^{t}}$   the equality (\ref{eee}). For each $t$ there appears
a summand
$$(-1)^{\frac{k-2}{2}-1}\check{z}_{n0}pPfF_{I^t\setminus
\{-n0\}}=(-1)^{\frac{k}{2}-2}\check{z}_{n0}pPfF_{I\setminus \pm
i,0,-n}.$$

Also there appear summands  $$ PfF_{I^t\setminus\{\pm i_s\}}=\pm
PfF_{I\setminus\{\pm i_t,\pm i_s\}}.$$  But the sum of these
summands over $t$ and $s$ is zero. Let $0<t<s$. If this summand
comes from the summand $ PfF_{I\setminus\{\pm i_s\}}$ in
(\ref{eee}), then it appears with the sign
$(-1)^{\frac{k}{2}-s-1}(-1)^{(\frac{k}{2}-1)-t-1}$.  If it comes
from the summand $ PfF_{I\setminus\{\pm i_t\}}$ in (\ref{eee})  then
it has the sign
$(-1)^{\frac{k}{2}-t-1}(-1)^{(\frac{k}{2}-1)-(s-1)-1}$. The sum of
these signs is zero.

Hence  in $Z(\mathfrak{o}_{2n+1},\mathfrak{o}_{2n-1})$ one has
$$PfF_{I}=(-1)^{\frac{k}{2}-1}\check{z}_{n0}pPfF_{I\setminus
\{-n,0\}}-4\sum_{t=1}^{\frac{k}{2}}(-1)^{\frac{k}{2}-t-1}(-1)^{\frac{k}{2}-2}\check{z}_{n0}pPfF_{I\setminus
\pm i,0,-n}.$$

Apply the obtained formulae  to $I=\hat{n}$. Recall that according
to Proposition \ref{tps} one has $pPfF_{\widehat{ 0, \pm
n}}=D_{n-1}(F_{11},...,F_{(n-1)(n-1)})$, and $pPfF_{\widehat{ \pm
i,0, \pm
n}}=D_{n-2}(F_{11},...,\widehat{F_{ii}},...,F_{(n-1)(n-1)})$. Thus
one proves Theorem.
\endproof




\subsection{ Action of pfaffians in  the Gelfand-Tsetlin-Molev base.}
\label{dpb1}

 The Gelfand-Tsetlin-Molev base of a  $\mathfrak{o}_{2n+1}$-representation is a base of a Gelfand-Tsetlin type, whose construction is  based
   on restrictions
 $\mathfrak{o}_{2n+1}\downarrow \mathfrak{o}_{2n-1}$, in contrast to the classical Gelfand-Tsetlin base,
 whose construction is based on restrictions
  $\mathfrak{o}_N\downarrow \mathfrak{o}_{N-1}$.

One can give the following formal definition. The
Gelfand-Tsetlin-Molev base is a base of a
$\mathfrak{o}_{2n+1}$-representation, for different
  $n$  the procedures of constructions of bases must be coherent  in the following sense. A base of a  $\mathfrak{o}_{2n+1}$-representation
  is a union a bases in  $\mathfrak{o}_{2n-1}$-representation,
  into which a
  $\mathfrak{o}_{2n+1}$-representation splits when one restricts $\mathfrak{o}_{2n+1}\downarrow \mathfrak{o}_{2n-1}$.

All notations below are taken from \cite{Molev}.

Let $V$ be given a  $\mathfrak{o}_{2n+1}$-representation with the
highest weight $(\lambda_1,...,\lambda_n)$, where

$$0\geq \lambda_1\geq...\geq\lambda_n.$$  Base vectors  are indexed by tableaus $\Lambda$ of
type

$\lambda_{n,1},\lambda_{n,2},...,\lambda_{n,n}$

$\sigma_n,\lambda'_{n,1},\lambda'_{n,2},...,\lambda'_{n,n}$

$\sigma_{n-1},\lambda_{n-1,1},\lambda_{n-1,2},...,\lambda_{n-1,n-1}$

...

$\lambda_{11}$

$\sigma_1,\lambda'_{11}$

The restrictions on these numbers are the following:
\begin{enumerate}

\item $\lambda_{ni}=\lambda_i$

\item $\sigma_k=0,1$

\item The inequalities hold:

$\lambda'_{k1}\geq\lambda_{k1}\geq\lambda'_{k2}\geq...\geq\lambda'_{k,k-1}\geq\lambda'_{kk}\geq\lambda_{kk}$
for $k=1,...,n$, and

$\lambda'_{k1}\geq\lambda_{k-1,1}\geq\lambda'_{k2}\geq...\geq\lambda'_{k,k-1}\geq\lambda_{k-1,k-1}\geq\lambda'_{kk}$
for $k=2,...,n$.

\item If $\lambda'_{k1}=0$, then $\sigma_k=0$.
\end{enumerate}

The first row is the highest weight of a
$\mathfrak{o}_{2n+1}$-representation $V$.

The third row is a weight of a  $\mathfrak{o}_{2n-1}$-representation
that appears if one restricts $\mathfrak{o}_{2n+1}\downarrow
\mathfrak{o}_{2n-1}$ and that contains the base vector.

The second row is set of indices of base vectors in $V_{\mu}^+$. And
so on.

Now we obtain formulae of the action of $PfF_{\widehat{n}}$ on a
base in $V_{\mu}^{+}$ and then in the Gelfand-Tsetlin-Molev base in
$V$.





Write the equality $V=\sum_{\mu}V_{\mu}^+\otimes V(\mu)$. From one
other hand, the pfaffian commutes with $\mathfrak{o}_{2n-1}$. Thus
the action on $V=\sum_{\mu}V_{\mu}^+\otimes V(\mu)$
 is written as
$\sum_{\mu}(PfF_{\widehat{n}}\mid_{V_{\mu}^+})\otimes id$.
 From the other using the theorem we relate the action of $PfF_{\widehat{n}}$
 on $V_{\mu}^+$ with the action of $z_{0n}$. The action of  $z_{0n}$
 on $V_{\mu}^+$  is described in
\cite{Molev}.

Using the description we obtain the following theorem.

Let $\rho_i=i-\frac{1}{2}$ for $i>0$, and $\rho_{-i}=-\rho_{i}$,
$$\gamma_i=\lambda'_{2i}+\rho_i+\frac{1}{2},$$

 $$\overline{\sigma}=\sigma+1 \,mod\, 2.$$

\begin{thm}
\label{maint} On the vector $\xi_{\Lambda}$ the pfaffian
$PfF_{\widehat{2}}$ acts as follows.

Let $\Lambda=(\lambda,\sigma,\lambda',\Lambda')$, where $\lambda$ is
the first row of   $\Lambda$, $\sigma,\lambda'$ is the second row of
$\Lambda$ and $\Lambda'$ is the remaining part of $\Lambda$.

If $\sigma=0$, then
$$PfF_{\widehat{2}}\xi_{\sigma,\lambda',\Lambda'}=C\xi_{\bar{\sigma},\lambda,\Lambda'},$$

where $C$ is some constant.

If $\sigma=1$, then
$$PfF_{\widehat{2}}\xi_{\sigma,\lambda,\Lambda'}=(-1)^n\sum_{j=1}^n\Pi_{t=1,t\neq
j}^n\frac{-\gamma_t^2}{\gamma_j^2-\gamma_t^2}\xi_{\bar{\sigma},\lambda'+\delta_j,\Lambda'},$$

Where $\lambda'+\delta_j$ is the row $\lambda'$ with $1$ added to
the $j$-th component.

Here
$C=\frac{C_n(\lambda_{n-1,1},...,\lambda_{n-1,n-1})}{\prod_{i=1}^{n-1}(\lambda_{n-1,i}+i-1)}$(see
the definition \ref{ch}).

\end{thm}

Note that in the case $n=2$, which corresponds to $\mathfrak{o}_5$
one has $C=\frac{C_2(\lambda_{1,1})}{(\lambda_{1,1})}=1$ and we
obtain the theorem \ref{maint1}.

\section{Conclusion}

We have constructed a quantum number $k$. Using it together with the
quantum numbers $N$ (the number of particles), $T,\tau_0$ (the
isospin and its projection) we can  classify the states of the
five-dimensional quasi-spin. The quantum number $k$ is defined
through its creation operator. As the creation operator we use the
noncommutative pfaffian associated with the algebra
$\mathfrak{o}_5$, which is isomorphic to the quasi-spin algebra. We
describe the action of the pfaffian on the other quantum numbers. In
particular the pfaffian increases by one the number of particles and
conserves the isospin and its projection.

\section{Acknowledgements}
 The work is supported by grants  Nsh-5998.2012.1,  MK-1378.2009.1.

\end{document}